\documentclass{article}
\usepackage{amsmath}
\usepackage[left=2.96cm, right=2.96cm, bottom=2.6cm, top=2.6cm]{geometry}
\usepackage{hyperref}
\usepackage{cleveref}
\usepackage{amssymb}
\usepackage{mathtools}
\usepackage{tikz-cd}
\usepackage{amsthm}
\usepackage{graphicx}
\usepackage{graphics}
\usepackage{mathrsfs}
\usepackage{cleveref}
\usepackage{thmtools}
\usepackage{enumitem}
\usepackage{bbm}
\newlist{steps}{enumerate}{1}
\usepackage{mathtools}
\makeatletter
\newcommand{\xMapsto}[2][]{\ext@arrow 0599{\Mapstofill@}{#1}{#2}}
\def\Mapstofill@{\arrowfill@{\Mapstochar\Relbar}\Relbar\Rightarrow}
\makeatother
\usepackage{array}
\setlist[steps, 1]{label = Step \arabic*:}
\theoremstyle{plain}
\usepackage{footnote}
\makesavenoteenv{tabular}
\newtheorem*{theorem*}{Theorem}
\newtheorem{theorem}{Theorem}[subsection]
\newtheorem{definition}[theorem]{Definition}

\newtheorem{remark}[theorem]{Remark}
\newtheorem{corollary}[theorem]{Corollary}
\newtheorem*{corollary*}{Corollary}
\newtheorem{exmp}[theorem]{Example}

\newtheorem{lemma}[theorem]{Lemma}
\newtheorem*{lemma*}{Lemma}

\numberwithin{equation}{subsection}
\usepackage{rotating}
\usepackage{comment}
\usepackage[nottoc,notlot,notlof]{tocbibind} 

\usepackage[toc,page]{appendix}

\title{Geometric Quantization Without Polarizations} 
\author{Joshua Lackman\footnote{josh@pku.edu.cn} }
\date{}
\begin{document}
\maketitle
\begin{abstract}
\noindent We derive the quantization map in geometric quantization of symplectic manifolds via the Poisson sigma model. This gives a polarization-free (path integral) definition of quantization which pieces together most known quantization schemes. We explain how this allows Schur's lemma to address the invariance of polarization problem. We compute this quantization map for the torus and obtain the noncommutative torus and its standard irreducible representation. 
\end{abstract}
\tableofcontents
\section{Introduction}
There are many different quantization schemes, including formal ones. To name some: there is formal deformation quantization of Fedosov \cite{fedosov} and Kontsevich \cite{kontsevich}, where the latter was explained in \cite{catt} using the Poisson sigma model; there is Kostant—Souriau's geometric quantization of symplectic manifolds \cite{woodhouse}; there is Weinstein's geometric quantization of Poisson manifolds via symplectic groupoids \cite{weinstein}, \cite{eli} (related are \cite{kar}, \cite{zak1}, \cite{zak2}).\footnote{There is also the A-model approach of Gukov—Witten \cite{guk}, which we discuss in \cref{A}.} The formal schemes, which produce star products, always work and are the only completed programs. However, they don't produce a $C^*$-algebra or states a priori, which is important physically and mathematically. In this paper, we show that the aforementioned quantization schemes all fit together (consistent with the concluding remarks of \cite{bon}). 
\\\\We will use path integrals of the Poisson sigma model to derive the quantization map, which sends functions on phase space to sections of the prequantum line bundle over the symplectic groupoid\footnote{Thus, observables are identified with states of another quantum theory, as in the A-model approach.} (this is the reduced phase space of the PSM),\footnote{For $(M,\omega)$ simply connected, the symplectic groupoid is $(M\times M,\pi_2^*\omega-\pi_1^*\omega).$} which act on (polarized) sections of the prequantum line bundle over the base — sections of the prequantum line bundle over this groupoid define operators on the prequantum Hilbert space, via their identification with integral kernels. A symplectic connection (used in Fedosov's approach) can be used to put these path integrals on a lattice, and one can try to take the limit as the spacing goes to zero (appendix \ref{appen}). In the more general Poisson case, we discussed this in \cite{Lackman3}.
\\\\In the geometric quantization of symplectic manifolds, there is an implicit \textit{assumption} that polarized sections form irreducible representations of some $C^*$-algebra. In principle, we don't need to take the representation to consist of polarized sections. However, the $C^*$-algebra in Kostant-Souriau's quantization scheme is \textit{defined} to be the operators on polarized sections (the groupoid approach \cite{weinstein}, \cite{eli} has a similar polarization dependence), and it quantizes only the very small subspace of functions which preserve polarized sections. In the case of $\textup{T}^*\mathbb{R}$ with its standard K\"{a}hler structure, \textit{only} affine functions preserve polarized sections for all complex affine polarizations.
\\\\The primary role given to Lagrangian polarizations is common to most quantization schemes,\footnote{The A-model approach is an exception. This approach also demotes the role of polarizations.} which introduces difficulties for several reasons. Aside from not providing a well-defined quantization map (or a well-defined classical limit as a result), they don't always exist (\cite{go}) and it can be difficult to relate quantizations obtained using different polarizations (\cite{brane}), which is largely caused by the $C^*$-algebra's dependence on the polarization. Additionally, when there is a Bohr—Sommerfeld condition one can't use continuous sections, and sometimes singular polarizations are used (\cite{bone}), bringing into question what a polarization is (see \cref{commu}). On the other hand, pure states (or irreducible representations) of $C^*$-algebras \textit{always} exist, which suggests it should be constructed more abstractly, as we aim to do.
\\\\Due to the difficulties and ambiguities introduced by treating polarizations as fundamental, we wish to demote their role and give an alternative perspective on quantization. The quantization should act on the entire prequantum Hilbert space, which acquires the physical interpretation of being a Hilbert space of mixed states (eg. density matrices define mixed states). Polarizations are just one tool for finding pure states: a polarization $\mathcal{P}$ which satisfies\footnote{This is similar to the use of charge in QFT for picking a superselection sector of the Hilbert space, see \cref{charge}.} 
\begin{equation}\label{commu}
[Q_f,\nabla_X]=0
\end{equation}
for all observables $Q_f$ and $X\in\mathcal{P}$ defines a subrepresentation of the quantization.\footnote{For any polarization, the Kostant-Souriau prequantization map only satisfies this for a small class of functions.}
This allows Schur's lemma (together with the BKS pairing) to help address the invariance of polarization problem. First, we should define exactly what a quantization map is — we take Rieffel's $C^*$-algebra perspective (\cite{rieffel}).\footnote{The abstract $C^*$-algebra perspective is close to physicists' abstract Hilbert space approach to quantum theory.} Put simply, it's the non-perturbative form of a star product, see \ref{q}. This is a mathematically motivated approach which doesn't rely on the ``wave functions depend on half of the variables" heuristic.\footnote{Defining the Hilbert space before the $C^*$-algebra is defining a dual space $V^*$ before $V.$}
\\\\To explain how this idea can help address invariance of polarization (which isn't the main point), consider the following corollary of Schur's lemma:
\begin{lemma*}
Two irreducible subrepresentations $\mathcal{H}_1,\mathcal{H}_2\subset \mathcal{H}$ of a $C^*$-algebra representation are either orthogonal subspaces or there is a unique projective equivalence $\mathcal{H}_1\to\mathcal{H}_2$ intertwining the representations.
\end{lemma*}
In QFT, this lemma implies that different vacuum states are orthogonal. In geometric quantization, since the BKS pairing between differently polarized sections is nontrivial, we can use it to argue the following:
\begin{corollary*}
Fix a quantization which acts on the prequantum Hilbert space. If $\mathcal{P}_1,\mathcal{P}_2\subset T_{\mathbb{C}}M$ are polarizations whose polarized sections form irreducible subrepresentations, then there is a unique projective equivalence intertwining them.\footnote{There is a subtlety since real-polarized sections need not be normalizable, which we will address.}
\end{corollary*}
Such an argument only makes sense if the quantization is defined on the entire prequantum Hilbert space. It explains why some subgroups of symplectomorphisms lift to unitary equivalences between differently polarized sections — as we will see in the example below, there is a quantization which acts on the prequantum Hilbert space of $\textup{T}^*\mathbb{R}$ and fixes subspaces of sections polarized with respect to complex affine polarizations, and indeed, the action of $\textup{Sp}(2,\mathbb{C})$ lifts to unitary equivalences. See the bottom of page 9 in \cite{brane} for a discussion about this. That the entire group of symplectomorphisms can't lift is related to the fact that, in the lattice definition of the path integrals, one can't choose a 2-cocycle  which is invariant under symplectomorphisms (step 4 in appendix \ref{appen}).
\subsection{Quantization of $\textup{T}^*\mathbb{R}$ Without Polarizations}\label{cot}
Our approach is justified by results of the Poisson sigma model and involves computing path integrals, thus it has its own difficulties (which we believe are surmountable), and it helps provide a theoretical framework for understanding quantization. As a proof-of-concept, we briefly review the main example in \cite{Lackman3} (which contains more detail). We will quantize the symplectic torus later, obtaining the noncommutative torus and its standard irreducible representation.
\\\\Consider the pointed K\"{a}hler manifold $(\textup{T}^*\mathbb{R},\omega,I,0)$ with the canonical symplectic form. The quantization of a smooth function $f$ is a section of a trivializable line bundle over the symplectic groupoid, and can be identified with
\begin{align}\label{tor}
     Q_f:\textup{T}^*\mathbb{R}\times \textup{T}^*\mathbb{R}\to\mathbb{C}\;,\,\; Q_f(u,v)=\frac{1}{2\pi\hbar}\int_{\textup{T}^*\mathbb{R}}f(z)e^{\frac{i}{\hbar}P_0(u,v,z)}\,\omega\;,
\end{align}
where $P_0(u,v,z)$ is the (signed) area of the geodesic polygon with vertices $0,u,v,z \in \textup{T}^*\mathbb{R}\,,$ and the integration is over $z.$\footnote{We emphasize square-integrable functions, but this formula makes sense more generally.} These quantized functions generate our $C^*$-algebra (when quantizing using symplectic groupoids with the standard prescription (\cite{weinstein}, \cite{eli}), a $C^*$-algebra is obtained only after choosing a polarization and there is no prescribed quantization map).
\\\\With respect to the the non-perturbative Moyal product $\star$  (\cite{eli}, \cite{Zachos}) and the twisted convolution algebra $*$ on the groupoid (appendix \ref{symp}), 
\begin{equation}
    Q_{f\star g}=Q_f\ast Q_g\;.
    \end{equation}
Therefore, $Q$ is compatible with a star product.\footnote{The Kostant-Souriau prequantization map is not compatible with a star product. We view this as more mathematically motivated than the common argument against it, ie. ``wave functions only depend on half of the variables".}
\\\\The line bundle over the symplectic groupoid acts on the the prequantum line bundle over $\textup{T}^*\mathbb{R},$ and this determines a representation of these quantized functions as integral operators on the prequantum Hilbert space $L^2(\textup{T}^*\mathbb{R}).$ The result is 
\begin{align}\label{defi}
    Q_f\Psi(u)=\frac{1}{(2\pi\hbar)^2}\int_{\textup{T}^*\mathbb{R}\times \textup{T}^*\mathbb{R}}f(v)\Psi(z)e^{\frac{i}{\hbar}\Omega(u,v,z)}\,\omega\boxtimes\omega\;,
\end{align}
where $\Psi\in L^2(\textup{T}^*\mathbb{R})$ and $\Omega(u,v,z)$ is the (signed) area of the geodesic triangle determined by $u,v,z;$ the integral is over $v,z,$ with respect to the product measure. We have that 
\begin{equation}
     (Q_f\ast Q_g)\Psi=Q_f(Q_g\Psi)\;,
\end{equation}
so that this forms a representation of our $C^*$-algebra.
\\\\Polarized sections form irreducible subrepresentations for \textit{all} covariantly constant Lagrangian polarizations. That is, polarized sections with respect to complex affine polarizations, and the canonical connection $(pdq-qdp)/2,$ are of the form 
\begin{equation}
\Psi(p,q)=e^{\frac{i}{2\hbar}(ap+bq)(cp+dq)}\psi(ap+bq)
\end{equation} 
for $ad-bc=1.$ For such polarizations \cref{commu} is satisfied,\footnote{Any operator commuting with the representation will give a subrepresentation.} and it turns out that
\begin{equation}\label{simp}
Q_f\Psi(u)=\frac{1}{2\pi\hbar}\int_{\textup{T}^*\mathbb{R}}f(u')\Psi(u'-u)e^{\frac{i}{\hbar}\Omega(0,u,u')}\,\omega
\end{equation}
(this is a polarized section). 
In particular, we get the position, momentum and Segal-Bargmann (holomorphic) representations (\cite{hall}). On the small subspace of functions for which the Kostant–Souriau prescription preserves polarized sections, our quantization maps agree. In addition, our prescription agrees with the Weyl quantization on sections polarized along the projection map, and is related to Wigner's phase space approach to quantization via the Riesz representation theorem.\footnote{This is one of the few quantization approaches not based on Hilbert spaces.} See \cite{Moyal}, \cite{Wigner}, and \cite{Curtright} for an exposition. 
\begin{remark}
In \cref{defi}, $\Omega$ is a 2-cocycle on the groupoid, and cocycles determined by different K\"{a}hler structures agree in groupoid cohomology by the van Est isomorphism theorem (\cite{Crainic}, \cite{Lackman4}). The resulting $C^*$-algebras are isomorphic, however their commutants are different (\cref{commuu}), ie. they preserve different polarized sections.
\end{remark}
We have glossed over the fact that real-polarized sections need not be normalizable. This is usually fixed by introducing half-forms (see \cite{bates} for an exposition). However, as we will discuss in \cref{bbks}, there is (at least nearly) a canonical inner product on real-polarized sections determined by the BKS pairing.
\\\\\textbf{Question 1:} It is commonly implied in quantization that only polarized sections define pure states. Is this true of the representation in \cref{defi}? In other words, if a section of the prequantum line bundle generates an irreducible subrepresentation, is it polarized? Can a pure state be localized around a point in phase space? Does a Lagrangian submanifold with a flat section determine a pure state?\footnote{Some of the questions we ask may have known answers.}
\\\\This question doesn't make sense to ask in the traditional approach. In fact, it would appear that the answer to the first is no, by considering a state of the form
\begin{equation}
    e^{-\frac{(p^2+q^2)}{4\hbar}}\psi_1(p+iq)+ e^{-\frac{(p^2+4q^2)}{8\hbar}}\psi_2(p+2iq)\;.
\end{equation}
Now, by the Stone–von-Neumann theorem there is a unique irreducible representation of this $C^*$-algebra. Therefore, all pure states can be obtained from sections of the prequantum line bundle. This raises another question:
\\\\\textbf{Question 2:} Which mixed states can be obtained from sections of the prequantum line bundle? Can density matrices (ie. convex combinations of pure states) be obtained?
\subsection{Definition of Quantization}\label{q}
According to the definition we use, a quantization of a Poisson manifold is an injective quantization map into a $C^*$-algebra whose perturbative expansion is a star product. It sits between Rieffel's strict quantization (\cite{rieffel}) and his strict deformation quantization.\footnote{Of course, eventually one needs to pick a state. We will discuss this later.} In the following, one should think of $M_{\hbar}$ as a $C^*$-algebra with a parameter $\hbar.$
\begin{definition}\label{nd}(see \cite{eli2})
Let $(M,\Pi)$ be a Poisson manifold and let $A\subset[0,1]$ be a set containing $0$ as an accumulation point. For each $\hbar\in A,$ let $M_{\hbar}$ be a unital $C^*$-algebra such that $M_{0}=L^{\infty}(M)$ and let 
\begin{equation}
    Q_{\hbar}:C_c^{\infty}(M)\to M_{\hbar}
\end{equation}
be injective and $^*$- linear\footnote{This compatibility with $^*$ can be relaxed.} such that its image generates $M_{\hbar}\,.$ Furthermore, assume that $Q_0$ is the inclusion map. We say that $Q_{\hbar}$ is a (non-perturbative) deformation quantization of $(M,\Pi)$ if there is a star product $\star_{\hbar}$ on $C^{\infty}(M)[[\hbar]]$ such that, for all $n\in\mathbb{N},$
\begin{equation}
    \frac{1}{\hbar^{n}}||Q_{\hbar}(f)Q_{\hbar}(g)-Q_{\hbar}(f\star_{\hbar}^n g)||_{\hbar}\xrightarrow[]{\hbar \to 0} 0\;,
\end{equation}
where $f\star_{\hbar}^n g$ is the component of the star product up to order $n.$
\end{definition}
If $Q_{\hbar}$ has an algebraically closed image, then we will get a genuine deformation of $C_c^{\infty}(M).$ However, we believe this happens only in very special cases, eg. in the symplectic case it seems to only happen for quotients of $\textup{T}^*\mathbb{R}^n$ by discrete groups, see \hyperref[q5]{question 5}. The compatibility with the star product means that the image of $Q_{\hbar}$ is \textit{nearly} closed.
\\\\In \cite{schl} (see also \cite{bord}) it is shown that there is a sequence of ``quantization maps" associated to a compact K\"{a}hler manifold (Berezin-Toeplitz quantization) which produces a formal deformation quantization in the sense of \ref{nd}, except that they have infinite dimensional kernel and finite dimensional image. We expect that these maps factor through an injective quantization map.
\\\\\textbf{Question 3:} Suppose we have a quantization map as in \cref{nd}. If we forget the quantization map but remember the family of $C^*$-algebras, can we recover the Poisson manifold? In other words, if a classical limit exists, is it unique?
\\\\Furthermore,
\\\\\textbf{Question 4:} Given a (faithful) representation of a quantization on the prequantum line bundle, does there always exist an irreducible subrepresentation?
\subsection{A Brief Review of $C^*$-Algebras and Invariance of Polarization}\label{review}
We briefly review the relevant theory of $C^*$-algebras and their states, including their relation to vector states, polarizations and Schur's lemma. See \cite{theo}, \cite{williams} for lecture notes on results in this section, \cite{murphy} for a textbook account and \cite{Gleason} for more on the relation with quantum mechanics.
\begin{definition}
A $C^*$-algebra $\mathcal{A}$ is a (complete) normed-algebra over $\mathbb{C},$ together with an anti-linear involution $^*,$ such that $\|AB\|\le\|A\|\|B\|$ and $\|A^*A\|=\|A^*\|\|A\|.$ A $C^*$-algebra is unital if it contains a unit, which we identify with $1.$
\end{definition}
\begin{definition}
A representation of a $C^*$-algebra $\mathcal{A}$ is a homomorphism into the bounded linear operators on a Hilbert space, denoted $\pi:\mathcal{A}\to\mathcal{B}(\mathcal{H}),$ for which $\pi(a^*)=\pi(a)^*$ for all $a\in\mathcal{A}.$
\end{definition}
\begin{definition}
A state on a $C^*$-algebra $\mathcal{A}$ is a linear functional $\rho:\mathcal{A}\to\mathbb{C}$ such that $\rho(1)=1$\footnote{We assume our $C^*$-algebras are unital.} and $0\le \rho(A^*A)\in\mathbb{R}.$
\end{definition}
Given a representation of $\mathcal{A}$ on a Hilbert space $\mathcal{H},$ a normalized vector $\Psi\in \mathcal{H}$ (which are called vector states) determines a state $\rho_{\Psi}$ by computing expectation values, ie.
\begin{equation}
    A\xmapsto{\rho_{\Psi}}\langle \Psi,A\Psi\rangle\;.
\end{equation}
Any state $\rho$ of $\mathcal{A}$ determines a vector state via the \textit{GNS construction}: consider the quotient $\mathcal{A}/\sim\,,$ where we quotient out by the left ideal consisting of those $A\in\mathcal{A}$ such that $\rho(A^*A)=0.$ There is a natural inner product on $\mathcal{A}/\sim\,,$ given by 
\begin{equation}
    \langle A,B\rangle=\rho(A^*B)\;.
    \end{equation}
This Hilbert space is a representation of $\mathcal{A},$ where $\langle 1,A1\rangle=\rho(A),$ so that the state determined by $1\in\mathcal{A}/\sim$ is equal to $\rho.$
\begin{definition}
A vector state of $\mathcal{A}$ (ie. a normalized vector in a representation of $\mathcal{A}$) is called pure if it generates an irreducible representation. Similarly, a state of $\mathcal{A}$ is pure if the associated vector state is pure.\footnote{This is really a theorem, but we'll use is as a definition. \textit{Pure} states are ones which can't be written as convex combinations of other states, thus are states of maximal information.}  Other states and vector states are called mixed.
\end{definition}
The fundamental states in quantum theory are pure states, whereas a density matrix determines a mixed state. If $\mathcal{H}$ is an irreducible representation of $\mathcal{A},$ then the map $\Psi\mapsto\rho_{\Psi}$ is an injection from rays in $\mathcal{H}$ into the space of pure states.
\\\\In classical theory (ie. commutative $C^*$-algebras of symplectic manifolds), a pure state is given by evaluation at a point, and the mixed states are normalized radon measures.
\\\\An elementary, but important result is that an intertwining map between irreducible representations of a $C^*$-algebra is either zero or is an isomorphism. This can be seen by noting that the kernel and image of an interwining map are subrepresentations.
\begin{definition}\label{commuu}
The commutant of a $C^*$-algebra representation $\pi:\mathcal{A}\to \mathcal{B}(\mathcal{H})$ consists of all $b\in \mathcal{B}(\mathcal{H})$ such that $\pi(a)b=b\pi(a),$ for all $a\in\mathcal{A}.$
\end{definition}
The following determines the commutant of an irreducible representation:
\begin{lemma}\label{schur}(Schur's lemma)
A representation of a $C^*$-algebra $\mathcal{A}$ on $\mathcal{H}$ is irreducible if and only if the commutant consists exactly of scalar multiples of the identity.
\end{lemma}
\Cref{schur} implies the following:
\begin{corollary}
Suppose that we have have irreducible representations of $\mathcal{A}$ on $\mathcal{H}_1,\mathcal{H}_2,$ such that there exist vector states $\Psi_1\in \mathcal{H}_1, \Psi_2\in \mathcal{H}_2$ which determine the same state of $\mathcal{A}.$ Then up to a constant in $S^1,$ there exists a unique unitary equivalence $\mathcal{H}_1\to\mathcal{H}_2$ intertwining the representations of $\mathcal{A}.$
\end{corollary}
Therefore, if there is any overlap in the states determined by two irreducible representations of a $C^*$-algebra, then the representations are projectively equivalent, in a unique way. The following is a corollary of Schur's lemma, it implies that there is an overlap of states if there is an overlap of vector states, with respect to the inner product. The map 
\begin{equation}
    T:\mathcal{H}_1\to\mathcal{H}_2
    \end{equation}
is the intertwining map obtained by restricting the inner product to $\mathcal{H}_2\times\mathcal{H}_1\to\mathbb{C}$ and using the Riesz representation theorem: 
\begin{corollary}\label{ortho}
If $\mathcal{H}$ is a representation of $\mathcal{A}$ and if $\mathcal{H}_1,\mathcal{H}_2$ are irreducible subrepresentations which are not orthogonal subspaces, then $T/\sqrt{T^*T}$ is a unitary equivalence, which induces the unique projective equivalence $\mathcal{H}_1\to\mathcal{H}_2$ intertwining the subrepresentations $(0<T^*T\in\mathbb{R}).$
\end{corollary}
Therefore, if our quantization is defined on all sections of the prequantum line bundle, and if we have two irreducible subrepresentatons $\mathcal{H}_1,\mathcal{H}_2$ that have nontrivial overlap, then there is a canonical unitary equivalence $\mathcal{H}_1\to\mathcal{H}_2.$ This suggests that there is a canonical unitary equivalence between differently polarized sections when they are irreducible subrepresentations of the same quantization,\footnote{The situation is a bit subtle for real-polarized sections, since technically, they may not live in the same Hilbert space. We discuss this more in \cref{ps}.} and it explains why the BKS pairing is sometimes a unitary equivalence (see \cref{bbks}).
\\\\Furthermore, the quantizations tend to surject onto the $C^*$-algebra of compact operators of polarized sections, which is, in fact, irreducible. In the simply connected case, the quantizations tend to be isomorphic to the $C^*$-algebra of compact operators of polarized sections, which is simple, meaning that there is a unique irreducible representation.
\begin{remark}\label{charge}
Since many vector states of the prequantum line bundle may determine the same state, one way of thinking of a polarization is as a gauge fixing by an intertwining operator. We will discuss polarizations more in \cref{ps}.
\end{remark}
\begin{remark}
It is interesting to observe that a measurement in quantum theory can naturally be modeled as the pulling back of a pure state of the quantization to the commutative $C^*$-subalgebra generated by an observable. A pure state pulls back to a (classical) mixed state, which is naturally identified with a probability measure via Gelfand duality. The same idea applies to any subset of operators which generates a commutative $C^*$-subalgebra.
\end{remark}
\section{The Quantization Scheme}
Here we will review the quantization scheme suggested for Poisson manifolds in \cite{Lackman3}, specialized to symplectic manifolds. First we motivate it (see appendix \ref{symp} for details about symplectic groupoids):
\\\\In \cite{kontsevich}, Kontsevich showed that all Poisson manifolds $(M,\Pi)$ admit a formal deformation quantization, ie. there is an associative (star) product $\star$ on $C^{\infty}(M)[[\hbar]]$ which deforms the product on $C^{\infty}(M),$ in such a way that 
\begin{equation}\label{star}
    f\star g=fg+i\hbar\{f,g\}+\mathcal{O}(\hbar^2)\;.
    \end{equation}
A Poisson sigma model approach to this star product was explained in \cite{catt}, which involves computing a 3-point function of operators inserted on the boundary of a disk. Specialized to the symplectic case $(M,\omega)$ and assuming that
\begin{equation}\label{condition}
    \int_{S^2}X^*\omega\in 2\pi\hbar\mathbb{Z}\;,
\end{equation}
it was shown in \cite{bon} (also argued in \cite{brane}, \cite{grady}) that 
\begin{equation}\label{fed}
    (f\star g)(m)=\int_{X(\infty)=m} f(X(1))g(X(0))\,e^{\frac{i}{h}S[X]}\,\mathcal{D}X\;,
\end{equation}
which is normalized so that $1\star 1=1.$ Here, the domain of integration is the space of contractible maps $X:S^1\to M,$\footnote{This should be definable non-perturbatively, it is essentially a 1-D sigma model.} where $0,1,\infty$ are three marked points, and the action is given by 
\begin{equation}
    S[X]=\int_D \tilde{X}^*\omega\;,
\end{equation}
where $\tilde{X}:D\to M$ is any disk agreeing with $X$ on the boundary. This is well-defined by \cref{condition}.
\begin{remark}\label{coc}
We can relate \cref{fed} to Fedosov's approach to formal deformation quantization: in order to define a star product in Fedosov's approach one needs to choose a symplectic connection. In order to define \ref{fed} on a lattice, one needs to choose a 2-cocycle on the local symplectic groupoid (ie. on a neighborhood of the diagonal in $M\times M$) which differentitates to $\omega$ under the van Est map, as explained in \cite{Lackman4}, see appendix \ref{appen}. The choice of a symplectic connection can be used to determine such a cocycle, and we get a lattice formulation of Fedosov's quantization.
\end{remark}
\textbf{Question 7:} Can one provide theoretical or numerical evidence that the lattice construction of \ref{fed} (given in appendix \ref{appen}) converges, on some manifolds? The same question goes for \cref{qant}.
\begin{definition}(formal)\label{formal}
We have an operator assignment \begin{equation}\label{ass}
   \mathcal{Q}:L^2(M,\omega)\to\mathcal{B}(L^2(M,\omega))\,,\footnote{This is the space of bounded linear operators on $L^2(M,\omega).$}\;f\mapsto\mathcal{Q}_f\,,\;\mathcal{Q}_fg=f\star g\;.
    \end{equation}
\end{definition}
This should give one way of defining a quantization map, assuming the lattice constructions converge as the spacing shrinks to zero. This quantization map is related to the quantization map we will describe on the symplectic groupoid.
\begin{remark}
There is little evidence that $\cref{fed}$ is associative non-perturbatively, and there is a no-go result on $S^2,$ see point 14 of \cite{rieffel} This is fine, as long as one interprets it as defining a $C^*$-algbera via the left action of $L^2(M,\omega)$ on itself, as in \cref{formal}. We don't expect to get an actual deformation of the product on $L^2(M,\omega).$ In other words, we view \ref{fed} as deforming the operator assignment $f\mapsto(g\mapsto fg),$ rather than as deforming the product $(f,g)\mapsto fg.$
\end{remark}
In light of this remark:
\\\\\textbf{Question 5:}\label{q5} Are there any examples (other than quotients of $\textup{T}^*\mathbb{R}^n$ by discrete groups) of symplectic manifolds for which one can deform the product (non-formally).
\\\\A negative result would show that \ref{fed} is not associative, non-perturbatively. It would also answer question nine of Rieffel's \cite{rieffel}, in the symplectic case.\footnote{This and similar questions are posed in \cite{cab}, which is somewhat related work.}
\subsection{Geometric Quantization of the Lie 2-Groupoid}
\Cref{fed} can be interpreted as a product in the twisted convolution of a Lie 2-groupoid, as described in \cite{Lackman3}. We need to do this in order to formulate our quantization scheme, which is easy to relate to Weinstein's. See appendix \ref{symp} for some related definitions.
\begin{definition}
We let $|\Delta^n|$ denote the geometric realization of the standard n-simplex, which can be identified with an $n$-ball with $(n+1)$-marked points on the boundary.
\end{definition}
We now define the 2-groupoid, which is a simplicial set that is a truncation of the singular simplicial set of $M$ in degree 2.
\begin{definition}(see \cite{zhuc})\label{lie2}
Let $(M,\omega)$ be a symplectic manifold. There is a Lie 2-groupoid, denoted $\Pi_{2}(M),$ which in degree $0$ is given by $M,$ in degree $1$ is given by
\begin{equation}
  \Pi_{2}^{(1)}(M)=\{\gamma\in \textup{Hom}([0,1],M\,): d\gamma\vert_{\{0,1\}}=0\}\;,
\end{equation}
and in degree $2$ is given by 
\begin{equation}
\Pi_{2}^{(2)}(M)=\{X\in \textup{Hom}(|\Delta^2|,M\,): dX\vert_{\textup{vertices}}=0\}/\sim \;,\footnote{$dX\vert_{\textup{vertices}}=0$ means the derivative vanishes on the vertices.}
\end{equation} 
where $\sim$ identifies two morphisms if they are homotopic relative to the boundary. We denote the face maps 
\begin{equation}
    \Pi_{2}^{(1)}(M)\to M
\end{equation}
by $s,t,$ where $s(\gamma)=\gamma(0),\,t(\gamma)=\gamma(1).$ We denote the
face maps 
\begin{equation}
   \Pi_{2}^{(2)}(M)\to \Pi_{2}^{(1)}(M)
\end{equation}
by $ d_0,\,d_1,\,d_2,$ where 
\begin{equation}
    d_0(X)=\gamma_1,\;d_1(X)=\gamma_2,\;d_2(X)=\gamma_0
    \end{equation}
are the 1-dimensional faces of $X.$ For such an $X,$ we may write 
\begin{equation}
    \gamma_0\cdot_X\gamma_1=\gamma_2.
    \end{equation}
\end{definition}$\,$
\\We can integrate the symplectic form $\omega$ to a 2-cocycle on the 2-groupoid, which defines a multiplicative line bundle $\mathcal{L}_2\to \Pi_{2}(M):$
\begin{definition}\label{multtt}
We have a 2-cocycle on $\Pi_2(M)$ given by 
\begin{equation}
    S:\Pi_{2}^{(2)}(M)\to\mathbb{R}\,,\;\; S[X]:=\int_{|\Delta^2|}X^*\omega\;.
\end{equation}
We then get a multiplicative line bundle 
\begin{equation}
    \mathcal{L}_2\to \Pi^{(1)}_{2}(M)\,,\;\;\mathcal{L}_2=\Pi_{2}^{(1)}(M)\times\mathbb{C}\;,
\end{equation}
whose multiplication is defined as follows: for $\gamma_0,\gamma_1$ composing to $\gamma_2$ via $X,$ we define
\begin{equation}
    (\gamma_0,\lambda_0)\cdot_X(\gamma_1,\lambda_1)=(\gamma_2,\lambda\beta \,e^{\frac{i}{\hbar}\int_{|\Delta^2|}X^*\omega}\,)\;.
\end{equation}
Sections of $\mathcal{L}_2\to  \Pi_{2}(M)$ are denoted $\Gamma(\mathcal{L}_2).$
\end{definition}
In the previous definition, if \cref{condition} holds, then the multiplication of the line bundle depends only on the boundary of $X.$ Now we define the twisted convolution algebra:
\begin{definition}
Consider two sections $w_1,\,w_2\in \Gamma(\mathcal{L}_2),$ ie. maps$\;\Pi_{2}^{(1)}(M)\to \mathbb{C}.$
Their twisted convolution
\begin{equation*}
w_1\ast w_2\in \Gamma(\mathcal{L}_2)
\end{equation*}
is given by\footnote{We are suppressing the dependence of $\ast$ on $\hbar$ in the notation.}
\begin{equation}
     (w_1\ast w_2)(\gamma)=\int_{\begin{subarray}{l}X\in \Pi^{(2)}_{2}(M)\\\gamma_2=\gamma\end{subarray}}w_1(\gamma_0)\cdot_X w_2(\gamma_1)\,\mathcal{D}X\;.
\end{equation}
Assuming \cref{condition}, we can take the integral to be over contractible maps $X:\partial|\Delta^2|\to M.$
 \end{definition}
Explicitly,
\begin{equation}
    (w_1\ast w_2)(\gamma)=\int_{\begin{subarray}{l}X\in \Pi^{(2)}_{2}(M)\\\gamma_2=\gamma\end{subarray}}w_1(\gamma_0)w_2(\gamma_1)e^{\frac{i}{\hbar}S[X]}\,\mathcal{D}X\;.
\end{equation}
This definition is formal and makes sense for any Lie 2-groupoid if we think of $X$ as just being a 2-morphism. If we apply it to a Lie 2-groupoid which is just a Lie 1-groupoid then we recover the usual definition of the twisted convolution algebra (see eg. \cite{eli}).
\subsection{Definition of the Quantization Map}
Here we define the quantization map. See appendix \ref{symp} for the relevant theory of symplectic groupoids and multiplicative line bundles, and appendix \ref{appen} for the lattice constructions of these path integrals.
\begin{definition}
We define the 2-groupoid quantization map by
\begin{equation}
    q:C^{\infty}(M)\to \Gamma(\mathcal{L}_2)\,,\;\;q_f(\gamma)=f(\gamma(1/2))\;.\footnote{See \cref{mult} for definition of $\mathcal{L}_2.$}
\end{equation}
\end{definition}
Looking at \cref{fed}, we immediately have:
\begin{lemma}
\begin{equation}
    f\star g=\iota^*(q_f\ast q_g)\;,
\end{equation}
where $\iota:M\to \Pi^{(1)}_2(M)$ is the map identifying an object with its associated identity morphism (ie. $m\in M$ is identified with the constant map $[0,1]\to M,\,t\mapsto m$).
\end{lemma} 
This is very close to Weinstein's geometric quantization of Poison manifolds, except that we are using a Lie 2-groupoid rather than the symplectic groupoid. However, we can push forward $q_f$ onto the symplectic groupoid. Before doing that, we recall its definition and quantization.
\\\\If the integrality condition \ref{condition} is satisfied, then the multiplicative line bundle of \cref{multtt} descends to a multiplicative line bundle over the source simply connected, symplectic groupoid $\Pi_1(M),$ which is the 1-truncation of $\Pi_2(M)$ (see appendix \ref{symp} for more). 
\begin{definition}\label{multi1}(\cite{ruif})
The arrows in $\Pi_1(M)\rightrightarrows M$ are paths $\gamma:[0,1]\to M$ up to homotopy, relative to the endpoints. The composition is given by concatenation. Letting $\widetilde{M}\to M$ be the universal cover of $M,$ we have 
\begin{equation}
    \Pi_1(M)\cong \widetilde{M}\times\widetilde{M}/\pi_1(M)\;,
\end{equation}
where we quotient out by the diagonal action of $\pi_1(M).$ The source and target are given by the projections. We denote such an arrow by $[\gamma].$ 
\end{definition}
The following is a standard construction in geometric quantization:
\begin{definition}
    Assuming the integrality condition \ref{condition}, the following defines a multiplicative line bundle with connection and Hermitian metric $(\mathcal{L}_1,\nabla,\langle \cdot,\cdot\rangle)\to \Pi_1(M)$ prequantizing $t^*\omega-s^*\omega:$
\\\\The vector space over an arrow $[\gamma]\in \Pi_1^{(1)}(M)$ 
is given by equivalence classes of points $(\gamma,\lambda)\in \mathcal{L}_2$ (defined in \cref{mult}). The equivalence  is such that if $\gamma_1$ is homotopic to  $\gamma_2$ relative to the endpoints, then
\begin{equation}
    (\gamma_1,\lambda)\sim (\gamma_2,\lambda\,e^{\frac{i}{\hbar}\int_D X^*\omega}\,)\;,
\end{equation}
where $X:D\to M$ is any homotopy between $\gamma,\gamma'$ and $\lambda\in\mathbb{C}.$ The multiplication of $\mathcal{L}_2$ descends to a multiplication of $\mathcal{L}_1.$ There is a canonical connection $\nabla$ on this line bundle, where if 
\begin{equation}
    X:[0,1]^2\to \Pi_1^{(1)}(M)
    \end{equation}
is a map agreeing with $\gamma_1,\,\gamma_2$ on $\{0\}\times[0,1],\,\{1\}\times[0,1],$ respectively, then parallel transport takes
\begin{equation}
    (\gamma_1,\lambda)\xmapsto{P^{\nabla}_{X}} (\gamma_2,\lambda\,e^{\frac{i}{\hbar}\int_{[0,1]^2}X^*\omega}\,)\;.
\end{equation}
\end{definition}
We now define the quantization map. In the following, $S^1$ has two marked points: the north and south pole. Therefore, $S^1$ is divided into two components.
\begin{definition}\label{qant}
Assuming the integrality condition \ref{condition}, we define the geometric quantization map by
\begin{equation}
    Q:L^2(M,\omega)\to \Gamma(\mathcal{L}_1)\;,\;\;Q_f([\gamma])=\big(\gamma,\,\int_{X} q_f(\gamma')e^{\frac{i}{\hbar}\int_D X^*\Pi}\,\mathcal{D}X\,\big)\;,
\end{equation}
where the integral is over contractible loops $X:S^1\to M$ with components $\gamma,\gamma'$ ($\gamma$ is fixed). This definition is independent of the representative $\gamma$ of $[\gamma].$
\end{definition}
This definition explains the presence of Fourier-like transforms in known cases where geometric quantization works, eg. Weyl quantization, \cite{bone}, \cite{weinstein}. In Weinstein's approach, the $C^*$-algebra is obtained by choosing a polarization on the symplectic groupoid and forming polarized sections. However, it doesn't prescribe a quantization map.
\begin{definition}
We can pull back the line bundle with connection and Hermitian metric 
\begin{equation}
    (\mathcal{L}_1,\nabla, \langle \cdot,\cdot\rangle)\to\Pi_1^{(1)}(M)
\end{equation} to the source fiber over $m\in M,$ which is the universal cover $\widetilde{M}$ of $M$ with basepoint $m.$ Denote this line  with connection and Hermitian metric by 
 \begin{equation}
    (\mathcal{L}_m,\nabla_m, \langle \cdot,\cdot\rangle_m)\to \widetilde{M}\;.
\end{equation}
There is a left representation of $\Gamma(\mathcal{L}_1)$ on itself. Therefore, there is an induced representation of $\Gamma(\mathcal{L}_1)$ on $\Gamma(\mathcal{L}_m),$ (see \cref{action}) and so our quantized functions act on this Hilbert space.
    \end{definition}
In the examples we know, when acting on $\Gamma(\mathcal{L}_m)$ it seems useful to divide $Q_f([\gamma])$ by 
\begin{equation}\label{kernel}
  \int_{m,\gamma(0),\gamma(1)} e^{\frac{i}{\hbar}S[X]}\,\mathcal{D}X\;,
\end{equation}
where the domain of integration consists of contractible maps $X:S^1\to M$ with the three marked points mapping to $m,\gamma(0),\gamma(1).$
\begin{remark}
If $M$ is prequantizable, then the representation of $Q_f$ on $\Gamma(\mathcal{L}_m)$ will descend to a representation on sections of the prequantum line bundle over $M.$ We will show this in the example of the torus. The examples we know are K\"{a}hler, for which there is a canonical symplectic connection, and therefore these path integrals should be canonically defined, see \cref{coc}. Sections which are polarized with respect to covariantly constant, complex polarizations seem to form (equivalent) irreducible subrepresentations, see \hyperref[q8]{question 8}.
\end{remark}
\begin{remark}
This quantization can be performed on (nice enough) Poisson manifolds (\cite{Lackman3}), in which case there is no ``measurement problem" between distinct isomorphism classes of objects in the 2-groupoid (ie. orbits/symplecitc leaves). This is because the quantizations of functions supported on disjoint orbits commute. For $\hbar=0,$ orbits are just points in phase space.
\end{remark}
\subsection{Remark: The 1-D Theory and Indefinite Pairing}
If $(M,\omega)$ is prequantizable to $(\mathcal{L},\nabla,\langle\cdot,\cdot\rangle)\to M,$ then we can write the action in \ref{fed} as the holonomy over the circle. Therefore, it is natural to consider maps with domain $[0,1],$ with the action given by parallel transport\footnote{This is related to the discussion starting on page 10 in \cite{brane}.} From this, we obtain an indefinite inner product on sections of the prequantum line bundle, given by
\begin{equation}\label{pair}
    \langle \Psi_0,\Psi_1\rangle=\int_{M\times M}\omega^n\boxtimes\omega^n\int_{X(0)=m_0}^{X(1)=m_1}\mathcal{D}X\,\big\langle\,e^{\frac{i}{\hbar}\int_{m_0}^{m_1}X^*\nabla}\Psi_0(m_0),\Psi_1(m_1)\rangle\;,
\end{equation}
where the path integral is over all maps $X:[0,1]\to M$ such that $X(0)=m_0, \,X(1)=m_1,$ and the exponential parallel transports $\Psi_0(m_0)$ to the vector space over $m_1.$ The kernel is
\begin{equation}\label{kernel}
    \int_{X(0)=m_0}^{X(1)=m_1}\mathcal{D}X\,e^{\frac{i}{\hbar}\int_{0}^{1}X^*\nabla}\;.
\end{equation}
In the case of $\textup{T}^*\mathbb{R},$ one could argue (see example \ref{ex}) that, up to a constant, this is equal to
\begin{equation}\label{b}
    e^{\frac{i}{2\hbar}(p_0q_1-q_0p_1)}\;,
\end{equation}
which is the parallel transport over the geodesic path from $m_0=(p_0,q_0)$ to $m_1=(p_1,q_1).$ \footnote{Associated to the K\"{a}hler structure is a symplectic connection, which determines a 1-cochain on the symplectic groupoid integrating $(p\,dq-q\,dp)/2.$ The result obtained using this 1-cochain to compute the path integral is \ref{b}, see appendix \ref{appen}.}
\\\\In quantum mechanics, $\mathcal{D}X\,e^{\frac{i}{\hbar}\int_{0}^{1}X^*\nabla}$ is normally integrated over paths with endpoints on Lagrangian submanifolds (with a Hamiltonian included), which reflects the fact that \ref{pair} seems to only be positive definite on pure states, in which case it agrees with the standard inner product on sections. In particular, it agrees with the standard inner product on K\"{a}hler polarized states. 
\\\\\textbf{Question 6:} Is a subrepresentation irreducible if and only if \cref{pair} pulls back to a positive definite form (and equals the inner product of sections of the prequantum line bundle)?
\section{Sympectic Torus}
We will describe the quantization map of the symplectic torus, which results in the noncommutative torus, and holomorphically polarized sections of the prequantum line bundle form a finite dimensional representation. The result is equivalent to working on the universal cover and using the example in \cref{cot} (see \cite{Lackman3} for more detail). See appendix \ref{symp} for some details on symplectic groupoids and their quantizations.
\\\\\textbf{Prequantization:}
\\The pointed K\"{a}hler structure $(\textup{T}^*\mathbb{R},n\omega,I,0)$ of \cref{cot} descends to a pointed  K\"{a}hler structure on the torus $T^2=\textup{T}^*\mathbb{R}/2\mathbb{Z},$ where we assume $n$ is a positive integer (which is only necessary to get a finite dimensional representation). In coordinates we have 
\begin{equation}
(p,q)\sim (p+2a,q+2b)\,,\;\;a,b\in\mathbb{Z}\;,
\end{equation}
with symplectic form and complex coordinate $dp\wedge dq,\;p+iq,$ and where $0$ is the point $(0,0).$
\\\\The  prequantization of $(\textup{T}^*\mathbb{R},\omega,I,0),$ consisting of the trivial line bundle with Hermitian metric and connection 
\begin{equation}
    \frac{p\,dq-q\,dp}{2}\;,
\end{equation} descend to a prequantization of $T^2.$ We can identify holomorphically polarized sections with smooth functions $\Psi:\mathbb{R}^2\to\mathbb{C}$ such that\footnote{These sections can be identified with theta functions (\cite{zas}) via a change of trivialization.}
\begin{align}\label{condi}
   \Psi(p+2a,q+2b)=e^{i\pi nab}e^{i\pi n(aq-bp)}\Psi(p,q)\,,\;\;\Psi(p,q)=e^{-\pi n(p^2+q^2)}\Psi_0(p+iq)\;.
\end{align}
where $\Psi_0$ is a holomorphic function on $\mathbb{R}^2.$ The theory of theta functions tells us that the dimension of the space of these sections is $n,$ which is the degree of the prequantum line bundle over $T^2.$
\\\\\textbf{Symplectic groupoid}
\\The (source simply connected) symplectic groupoid of $T^2$ is the action groupoid $\mathbb{R}^2\ltimes T^2\rightrightarrows T^2,$ where the source, target and multiplication are given by
\begin{align}
   &\nonumber s(x,y,p,q)=(p,q)\,,\;\;t(x,y,p,q)=(p+2x,q+2y)\,,
   \\& (x,y,p,q)\cdot(x',y',p+2x,q+2y)=(x+x',y+y',p,q)\;.
\end{align}
It's geometric quantization can be obtained using \cref{action}.
\\\\\textbf{Quantization Map}
\\Computing the $\star$ product \ref{fed} gives (see appendix \ref{appen}), for $f,g:T^2\to\mathbb{C},$
\begin{equation}\label{mt}
  (f\star g)(p,q)=  \frac{n^2}{4}\int_{\mathbb{R}^4}f(p_1,q_1)g(p_2,q_2)e^{i\pi n[(p_2-p)(q_1-q)-(q_2-q)(p_1-p)]}\,dp_1\,dq_1\,dp_2\,dq_2\;.
\end{equation}
One can check that this is well-defined and is just the Moyal product on the universal cover. We have that
\begin{align}
 & \nonumber e^{i\pi p}\star e^{i\pi q}= e^{\frac{i\pi}{n}}e^{i\pi p}e^{i\pi q}
    \\& e^{i\pi q}\star e^{i\pi p}= e^{-i\frac{\pi}{n}}e^{i\pi p}e^{i\pi q}\;,
\end{align}
Therefore, 
\begin{equation}
    e^{i\pi p}\star e^{i\pi q}= e^{\frac{2\pi i}{n}}e^{i\pi q}\star e^{i\pi p}
\end{equation}
and we have recovered the noncommutative torus (see \cite{weinstein} for the derivation using polarizations). Furthermore, $Q_f$ of \cref{tor} and its representation descends to a representation on holomorphically polarized sections of the prequantum line bundle.\footnote{We are using the fact that holomorphically polarized sections of the prequantum line bundle can be identified with sections over a source fiber of the symplectic groupoid of $T^2,$ satisfying \cref{condi}.} That is, for $f:T^2\to\mathbb{C},$ computing the representation of \cref{qant} gives
\begin{align}
   &\nonumber  Q_f\Psi(p,q)
   \\& =\frac{n^2}{4}\int_{\textup{T}^*\mathbb{R}\times \textup{T}^*\mathbb{R}}f(p_1,q_1)\Psi(p_2,q_2)e^{i\pi n[(p_2-p)(q_1-q)-(q_2-q)(p_1-p)]}\,dp_1\,dq_1\,dp_2\,dq_2\;,
\end{align}
which may be further simplified using \cref{simp}. That is, $Q_f\Psi$ satisfies \cref{condi} and furthermore  $Q_{f\star g}=Q_f\ast Q_g\,,$  $(Q_{f}\ast Q_{g})\Psi=Q_f(Q_g\Psi)\,,$ so that $Q$ is a quantization according to \cref{nd}. 
\section{Polarizations, States and the Classical Limit}\label{ps}
The quantization map has been defined independently of a polarization. The only data that should be needed to define it is a 2-cocycle on the local symplectic groupoid which differentiaties to the symplectic form under the van Est map (see appendix \ref{appen}). Such a cocycle always exists and can be determined by a symplectic connection (a torsion-free connection for which the symplectic form is constant). In addition, the symplectomorphisms which lift to the quantization are going to be limited by the 2-cocycle. However, different 2-cocycles may lead to isomorphic quantizations. 
\\\\Once a quantization has been computed, one needs to choose a Hamiltonian and a state in order to do quantum mechanics. The time evolution can be computed independently of any state, using the Heisenberg picture. A mixed state can be used if one wants to describe a statistical ensemble, otherwise a pure state should be used – these are states which generate irreducible representations via the GNS construction, see \cref{review}. Pure states always exist, so one doesn't need need to commit to polarized sections, which may not exist. In the case of Poisson manifolds, \textit{many} distinct irreducible representations can be determined by different subspaces of sections of the line bundle.
\\\\When they exist, polarized sections do seem to often form subrepresentations of the $C^*$-algebra determined by our quantization map, but it would be nice to have a precise statement. In the examples we know, for each Lagrangian submanifold in the polarization, the 2-cocycle vanishes on composable pairs of arrows which each have source and target on the Lagrangian submanifold. This is a local analogue of the infinitesimal Lagrangian condition.
\\\\\textbf{Question 8:}\label{q8} What is the compatibility condition needed between the 2-cocycle and the polarization to ensure that polarized sections form subrepresentations? Is it enough that the cocycle is determined by a symplectic connection and that it satisfies the local Lagrangian condition discussed above? (for a discussion in the formal case, see \cite{nolle}, \cite{nolle2}).
\\\\There is no right quantization map\footnote{Though on a K\"{a}hler manifold there should be a canonical one.}, as there is no one-to-one correspondence between classical and quantum systems. Relatedly, there is a notion of gauge equivalence of star products, and every Poisson manifold has a canonical class. What quantization maps are good at is defining the classical limit of the quantum theory, as they endow the family of $C^*$-algebras with something like a topology. Let 
\begin{equation}
Q_{\hbar}:C_c^{\infty}(M)\to M_{\hbar}
\end{equation}
be our quantization, as in \cref{nd}. Consider a family of states $\rho_{\hbar}:M_{\hbar}\to \mathbb{C},$ ie. linear maps such that $\rho_{\hbar}(1)=1$ and $\rho_{\hbar}(H_{\hbar}^*H_{\hbar})\ge 0\in\mathbb{R},$ for all $H_{\hbar}\in M_{\hbar}.$ We can say that:
\begin{definition}
A sequence of observables $H_{\hbar}\xrightarrow[]{\hbar\to 0}H$ if $\|Q_{\hbar}(H)-H_{\hbar}\|_{\hbar}\to 0.$
\\\\Dually, $\rho_{\hbar}\xrightarrow[]{\hbar\to 0}\rho_0$ if for all $H\in C^{\infty}_c(M),$
\begin{equation}
    \rho_{\hbar}(Q_{\hbar}(H))\xrightarrow[]{\hbar\to 0}\rho_0(H)\;.
\end{equation}
\end{definition}
Classically, pure states are given by evaluation at a point and mixed states are given by integration with respect to normalized radon measures. The pure states determined by coherent states converge to a classical pure state.
\\\\\textbf{Question 9:} Which classical mixed states appear as limits of mixed vector states in the $C^*$-algebra of $\textup{T}^*\mathbb{R},$ with the representation on the prequantum line bundle?
\\\\Similarly, we can ask:
\\\\\textbf{Question 10:}
Suppose we have a quantization of a classical system. Do all classical states appear as limits of quantum states?
\subsection{States Determined by the BKS Pairing (Inner Products)}\label{bbks}
Assume that $(M,\omega)$ is symplectic with $\textup{dim} \,M=2n.$ If $M$ is noncompact then real-polarized sections need not be in the prequantum Hilbert space, and thus don't exactly inherit the inner product. However, assuming they form an irreducible subrepresentation of the underlying algebra of the $C^*$-algebra, there should, at least \textit{nearly}, be a canonical inner product. We can explain this as follows: \\\\Let $\Psi_1$ be a section of the prequantum line bundle which generates an irreducible representation of the underlying algebra of the $C^*$-algebra. Suppose that one can find a second section $\Psi_2$ which generates an irreducible subrepresentation, such that:
\begin{equation}\label{bks}
    \int_{M}\langle \Psi_2,\Psi_1\rangle\, \omega^n=1
\end{equation}
and
\begin{equation}\label{2}
   \int_{M}\langle \Psi_2,A^*A\Psi_1\rangle\, \omega^n\ge 0\in\mathbb{R}
\end{equation}
for all operators $A$ in our $C^*$-algebra. Then the pair $\Psi_1,\Psi_2$ define a pure state. The pairing \ref{bks} is the BKS pairing\footnote{See \cite{gui}, and \cite{bates} page 108 for an exposition.} if $\Psi_1,\Psi_2$ are polarized. Using the GNS construction (see \cref{review}), we get an inner product on the sections generated by $\Psi_1,$ defined by
\begin{equation}
    \langle A\Psi_1, B\Psi_1\rangle =\int_{M}\langle \Psi_2,A^*B\Psi_1\rangle\, \omega^n\;,
\end{equation}
and this is a $C^*$-algebra representation.
\begin{exmp}
Consider the example of $\textup{T}^*\mathbb{R},$ and let $\Psi_1(p,q)e^{-\frac{i}{2\hbar}pq}\Psi(q),$ ie. a position space wave function. Let 
\begin{equation}
    \Psi_2(p,q)=e^{\frac{i}{2\hbar}pq}\int_{-\infty}^{\infty}\Psi(q)\,e^{-\frac{i}{\hbar}pq}\,dq\;.
    \end{equation}
Then up to normalization, the pair $\Psi_1, \Psi_2$ satisfies \ref{bks}, \ref{2}. The induced inner product is the pairing with respect to the Lebesgue measure.
\end{exmp}
We can make a stronger statement if we assume a complex polarization exists for which the polarized sections form an irreducible subrepresentation: suppose that $\{\Psi_{1_i}\}_i$ are (real or complex) polarized sections and $\{\Psi_{2_i}\}_i$ are complex-polarized sections, both which form irreducible subrepresentations. The sections $\{\Psi_{2_i}\}_i$ belong to the prequantum Hilbert space and the BKS pairing, together with the Riesz representation theorem, define a linear map 
\begin{equation}
    \{\Psi_{1_i}\}_i\to \{\Psi_{2_i}\}_i\;.
\end{equation} 
For each $\Psi_{1_i}$ this map determines a $\Psi_{2_j}$ which satisfies \cref{bks}, \ref{2}. The resulting Hilbert space structure on $\{\Psi_{1_i}\}_i$ is unitarily equivalent to $\{\Psi_{2_i}\}_i$ and is independent of the complex polarization by \cref{ortho}.
\begin{remark}
Suppose we have an irreducible representation of a $C^*$-algebra on $(\mathcal{H},\langle\cdot,\cdot\rangle_1).$ Suppose that $\langle\cdot,\cdot\rangle_2$ is another inner product on $\mathcal{H}$ whose norm is equivalent to $\langle\cdot,\cdot\rangle_1$ and which is also compatible with the involution $^*.$ Then by Schur's lemma, $\langle\cdot,\cdot\rangle_2=\lambda\langle\cdot,\cdot\rangle_1$ for some $\lambda>0.$ So if we could find one good inner product on real-polarized sections, it is probably the one we want.
\end{remark}
We can formalize this construction by considering a $C^*$-algebra together with a representation of the underlying algebra on a vector space $V,$ equipped with a partially defined inner product (\cite{ant}) on a (symmetric) subspace of $V\otimes V$ which is compatible with the involution. Furthermore, we require that vectors with finite norm form a Hilbert space $\mathcal{H}$ and that the inner product is defined on $V\otimes\mathcal{H},$ and that it varies continuously in $\mathcal{H}.$
\section{Comparison With the A-Model}\label{A}
Here, we briefly survey some known connections with the A-model. In \cite{catt1}, it is explained that a certain gauge fixing of the Poisson sigma model on closed surfaces results in the action of the A-model (after a partial integration).
\\\\In the case of the disk, both our approach and the A-model\footnote{The A-model approach is for non-degenerate Poisson manifolds, ie. symplectic, which have a complexification $Y$ with a good A-model, eg. $Y$ is hyper-K\"{a}hler. It aims to quantize only polynomials. } provide frameworks for understanding quantization which demote the role of Lagrangian polarizations, with some of the same goals. In the A-model approach \cite{guk}, \cite{brane}, one first complexifies the symplectic manifold $(M,\omega_M)$ to $(Y,\omega_Y)$ (see \cite{aldi}, \cite{francis} for related work\footnote{In \cite{aldi}, they also obtain a representation of the noncommutative torus, using $M=T^2, Y=T^4.$ In \cite{francis},  they quantize toric Poisson varieties, where $Y$ is given by the symplectic groupoid and $M$ is the (Lagrangian) identity bisection.}). \textit{Roughly}, the observables are given by $\textup{Hom}(\mathcal{B}_{cc},\mathcal{B}_{cc}),$ where $\mathcal{B}_{cc}$ is the canonical coisotropic brane supported on $Y,$ and the multiplication is given by integrating the integrand of the path integral $\ref{fed}$ over maps 
\begin{equation}
X:S^1\to Y
\end{equation}
which bound holomorphic disks (page 21 of \cite{brane}). On the other hand, we've been integrating it over maps 
\begin{equation}
    X:S^1\to M\xhookrightarrow{} Y
    \end{equation}
which bound disks. Thus, the difference between the two is which cycle in the space of maps $X:S^1\to Y$ is being used in the domain of the path intergal. Their Hilbert space is given by $\textup{Hom}(\mathcal{B},\mathcal{B}_{cc}),$ where $\mathcal{B}$ is a flat line bundle over $M.$ In some cases this Hilbert space is supposed to be identified with a Hilbert space of polarized sections, which means that we have two quantizations acting on it. It would be interesting to know what the relationship is.\footnote{The observables in this approach are states of the A-model, similarly to how in the approach we use the observables are states in the Hilbert space of the symplectic groupoid.}
\\\\Additionally, there is an idea that the Fukaya category of a symplectic groupoid should be monoidal and act on the Fukaya category of the symplectic leaves of the base manifold (\cite{Pa}). This is closely related to the fact that observables of the Poisson manifold should (in nice enough cases) define states in the quantization of the symplectic groupoid, which act on the quantizations of the symplectic leaves (\cite{Lackman3}, \cite{wein}).
\\\\Lastly, one should note that Fukaya used the holonomy groupoid in his approach to the $A_{\infty}$-category of Lagrangian polarizations \cite{fuk}. On $T^2,$ the composition of this category is formally the same as the (non-perturbative) Moyal product \cref{mt} (an explanation of this is asked in question 6 of \cite{soi}). Fukaya indicates that the relationship between the $C^*$-algebra of a polarization and the $C^*$-algebra of a quantization is an instance of a mirror symmetry\footnote{In homological mirror symmetry, Lagrangian polarizations are also given a primary role (\cite{deni}).} (page 8 of \cite{fuk}).
\begin{appendices}
\section{Symplectic Groupoids and Multiplicative Line Bundles}\label{symp}
Because we are only working with Poisson manifolds which are symplectic in this paper, thee relevant theory of symplectic groupoids simplifies. We will briefly go over the details we need. See \cite{eli} for more, as well as \cite{ruif}.
\begin{definition}
A groupoid is a category $G\rightrightarrows X$ for which the objects $X$ and arrows $G$ are sets and for which every morphism is invertible. Notationally, we have two sets $X, G$ with structure maps of the following form:
\begin{align*}
    & s,t:G\to X\,,
   \\ & \iota:X\to G\,,
    \\ & m:G\sideset{_t}{_{s}}{\mathop{\times}} G\to G\,,
    \\& ^{-1}:G\to G\,.
\end{align*}
Here $s,t$ are the source and target maps; $\iota$ is the identity bisection (ie. $X$ can be thought of as the set of identity arrows inside $G$); $m$ is the multiplication, which for short may be denoted $(g_1,g_2)\mapsto g_1\cdot g_2,$ and $^{-1}$ is the inversion map. We frequently identify a point $x\in X$ with its image in $G$ under $\iota$ and write $x\in G.$  Sometimes the space of arrows is denoted $G^{(1)},$ with $G^{(2)}$ the space of composable arrows. 
\\\\A Lie groupoid is a groupoid $G\rightrightarrows X$ such that $G, X$ are smooth manifolds, such that all structure maps are smooth and such that the source and target maps submersions.
\end{definition}
\begin{definition}
A symplectic groupoid is a Lie groupoid with a symplectic form $\omega$ on its space of arrows, such that $\pi_1^*\omega+\pi_2^*\omega-m^*\omega=0,$ where $\pi_1,\pi_2:G\sideset{_t}{_{s}}{\mathop{\times}} G\to G$ are the projections.
\end{definition}
The manifold of objects of a symplectic groupoid is canonically a Poisson manifold, and the groupoid integrates the Lie algebroid associated to the Poisson manifold.
\begin{definition}\label{fundg}
Let $(M,\omega)$ be a symplectic manifold. The fundamental groupoid, denoted 
\begin{equation}
    \Pi_1(M)\rightrightarrows M\;,
    \end{equation}
is the symplectic groupoid whose arrows between two objects $m_1,m_2$ are homotopy classes of smooth maps starting at $m_1$ and ending at $m_2,$ with multiplication given by concatenation. The space of arrows is naturally identified with 
\begin{equation}
    \widetilde{M}\times\widetilde{M}/\pi_1(M)\;,
    \end{equation}
where $\widetilde{M}$ is the universal cover of $M$ and $\pi_1(M)$ acts by the diagonal action. The identity objects are identified with constant paths and the symplectic form is given by $t^*\omega-s^*\omega.$
\\\\If $M$ is simply connected, then this groupoid is called the pair groupoid, denoted $\textup{Pair}\,M\rightrightarrows M.$ This is the simplest kind of groupoid and its space of arrows is given by $M\times M,$ with the source and target given by the projections and multiplication given by $(m_1,m_2)\cdot(m_2,m_3)=(m_1,m_3).$ The identity bisection $\iota$ is the embedding into the diagonal.
\end{definition}
\begin{definition}\label{action}
If $(M,\omega)$ is prequantizable with line bundle $\mathcal{L}\to M,$ then 
\begin{equation}\label{mult}
    t^*\mathcal{L}\otimes s^*\mathcal{L}^*\to \Pi_1(M)
    \end{equation}
prequantizes the space of arrows of its fundamental groupoid. The vector space over $[\gamma]\in \Pi_1(M)$ is naturally identified with the space of linear maps between the vector spaces over its source and target in $\mathcal{L}.$ Therefore, there is a natural multiplication of vectors over a composable pair of arrows, given by composition of the linear maps. This turns \ref{mult} into a multiplicative line bundle (ie. a line bundle over the space of arrows with a compatible multiplication). 
\\\\Due to the multiplication of the line bundle, sections of \ref{mult} have an associative product, called the twisted convolution. Assume the dimension of $M$ is $2n.$ For $w_1,w_2\in\Gamma(t^*\mathcal{L}\otimes s^*\mathcal{L}^*),$ it is given by
\begin{equation}
(w_1\ast w_2)(\gamma)=\int_{\gamma'(1)=\gamma(1)} w_1(\gamma\cdot\gamma'^{-1})\cdot w_2(\gamma')\,t^*\omega^n\;.
\end{equation}
\\\\Now assume $M$ is simply connected. The sections of \ref{mult} act on sections of $\mathcal{L}$ via twisted convolution. Given $w\in\Gamma(t^*\mathcal{L}\otimes s^*\mathcal{L}^*),\,\Psi\in\Gamma(L),$ the action is given by
\begin{equation}
    \Psi\mapsto w\Psi\;,\;\;(w\Psi)(m)=\int_{M} w(m',m)(\Psi(m')) \,\omega^n\;.
\end{equation}
In the integrand we are using the identification of a vector in the multiplicative line bundle with a linear map between vector spaces in the line bundle over the base. This representation satisfies
\begin{equation}
    (w_1\ast w_2)\Psi=w_1(w_2\Psi)\;.
\end{equation}
\end{definition}
$\,$\\Note that, this multiplicative line bundle can exist even if $M$ isn't prequantizable, ie. it exists as long as $\int_{S^2}X^*\omega\in\mathbb{Z}$ for maps $X:S^2\to M.$
\section{Coordinate-Free Definition of Path Integrals on a Lattice}\label{appen}
Lie algebroids provide a natural setting for which path integrals can be put on a lattice. We will briefly describe the lattice construction given in \cite{Lackman4}, specialized to the example of the star product in the form \cref{fed}, ie.
\begin{equation}\label{sympp}
    (f\star g)(m)=\int_{X(\infty)=m} f(X(1))g(X(0))\,e^{\frac{i}{h}\int_{D}X^*\omega}\,\mathcal{D}X\;.
\end{equation}
This path integral is over the space of maps $\{X:D\to M: X(\infty)=m\},$\footnote{Maps $X:D\to M$ are equivalent to Lie algebroid morphisms $TD\to TM.$} where $D$ is the disk with marked points $0,1,\infty$ on the boundary.
The steps are:
\begin{enumerate}
    \item triangulate the disk $D$ (with the marked points contained in the vertices) and form the associated simplicial set $\Delta_D,$
    \item integrate the Poisson manifold to the local symplectic groupoid (ie. choose a small neighborhood of the diagonal in $M\times M$, which may depend on $\Delta_D$),
    \item approximate the space of continuous maps $D\to M$ by morphisms of simplicial sets between $\Delta_D$ and the local symplectic groupoid,
    \item integrate the symplectic form to a 2-cocycle on the local groupoid via the van Est map, $VE,$
    \item form the approximations of the action (generalized Riemann sums\footnote{This is a coordinate-free notion of Riemann sums which makes sense on manifolds (globally), and is needed to formalize lattice constructions of path integrals in \cite{Lackman4}, eg. Feynman's construction.}),
    \item construct the measure on the hom space of simplicial sets by using the symplectic form,
    \item compute the approximations of the functional integral,
    \item take a limit over triangulations of the disk (shrinking the neighborhoods in 2 as needed).
\end{enumerate}   
$\,$\\We can replace the disk with any $n$-manifold and the symplectic form with any $n$-form (more general data can also be used). Since $d(pdq-qdp)/2=dp\wedge dq,$ the following example applies to \ref{sympp} with symplectic $\textup{T}^*\mathbb{R}$ via Stokes' theorem, and consequently the torus:
\begin{exmp}\label{ex}
Consider the path integral in \cref{kernel} on $\textup{T}^*\mathbb{R}=\mathbb{R}^2,$ whose domain is $[0,1]:$
\begin{equation}\label{bc}
 \int_{X(0)=(p,q)}^{X(1)=(p',q')}\mathcal{D}X\,e^{\frac{i}{\hbar}\int_{0}^1X^*(pdq-qdp)/2}\;.   
\end{equation}
The symplectic groupoid is the pair groupoid $\textup{Pair}\,\mathbb{R}^2=\mathbb{R}^2\times\mathbb{R}^2\rightrightarrows\mathbb{R}^2.$ A $1$-cochain on the symplectic groupoid integrating $(pdq-qdp)/2$ is given by\footnote{This is step 4 in the lattice construction, ie. $VE(\Omega)=(pdq-qdp)/2.$} \begin{equation}
    \Omega(p_0,q_0,p_1,q_1)=(p_0q_1-q_1p_0)/2\;,
    \end{equation}
ie. differentiating $\Omega$ along $(p_1,q_1)$ at the diagonal gives $(pdq-qdp)/2.$ This is the 1-cochain obtained using the symplectic connection associated with the standard K\"{a}hler structure, ie. obtained by integrating $(pdq-qdp)/2$ over the convex hull of $(p_0,q_0), (p_1,q_1).$
\\\\A triangulation $\Delta_{[0,1]}$ of $[0,1]$ is determined by a sequence of points $0=t_0<t_1<\ldots<t_n=1.$ A morphism of simplicial sets
\begin{equation}
    X:\Delta_{[0,1]}\to\textup{Pair}\,\mathbb{R}^2
    \end{equation}
with the boundary conditions in \ref{bc} is just a sequence of points 
\begin{equation}
    (p,q),(p_1,q_1),\ldots,(p_{n-1},q_{n-1}),(p',q')\;,
    \end{equation}
and the measure on the space of morphisms is the product measure. The generalized Riemann sum associated to $X^*\Omega$ is
\begin{equation}
    \sum_{\Delta\in\Delta_{[0,1]}}X^*\Omega(\Delta)=\frac{1}{2}\sum_{i=0}^{n-1}(p_iq_{i+1}-q_{i+1}p_i)\;.
\end{equation}
Heuristicaly, we can compute \ref{bc} using the entire pair groupoid (as a physicist would do). Doing this using triangulations with an even number of vertices (and using the integral representation of the Dirac delta), up to a constant we get
\begin{equation}
     e^{\frac{i}{2\hbar}(p_0q_1-q_0p_1)}\;.
\end{equation}
\end{exmp}
\end{appendices}

 \end{document}